\documentclass[12pt]{article}
\usepackage{amssymb}
\usepackage{amssymb}
\usepackage{url}
\usepackage{amsmath, amscd}
\usepackage{latexsym}
\usepackage[dvips]{graphicx}
\usepackage{epstopdf}
\usepackage{lineno}
\usepackage{color}
\usepackage{setspace}
\usepackage{tikz}
\usepackage{tkz-graph}
\usepackage{tkz-berge}
\usepackage{multicol}
\usepackage{float}
\newtheorem{theorem}{Theorem}[section]
\newtheorem{lemma}[theorem]{Lemma}
\newtheorem{corollary}[theorem]{Corollary}

\newtheorem{proposition}[theorem]{Proposition}

\newcommand{\qed}{\hfill $\Box$ }
\newcommand{\proof}{\noindent{\bf Proof.}\ \ }
\begin{document}

\title{\Large {\bf Wiener indices of maximal $k$-degenerate graphs}}
\author{Allan Bickle\\
\normalsize  Department of Mathematics,\\
\normalsize  Penn State University, Altoona Campus,
\normalsize  Altoona, PA 16601, U.S.A.\\
\small {\tt E-mail: aub742@psu.edu}\\
\and Zhongyuan Che\\
\normalsize  Department of Mathematics,\\
\normalsize  Penn State University, Beaver Campus,
\normalsize  Monaca, PA 15061, U.S.A.\\
\small {\tt E-mail: zxc10@psu.edu}\\
}

\date{\today}
\maketitle

\begin{abstract} 
A graph is \emph{maximal $k$-degenerate} if each induced subgraph has a vertex of degree at most $k$ and adding any new edge 
to the graph violates this condition.
In this paper, we provide sharp lower and upper bounds on Wiener indices of maximal $k$-degenerate graphs of order $n \ge k \ge 1$. 
A graph is \emph{chordal} if every induced cycle in the graph is a triangle and chordal maximal $k$-degenerate graphs of order $n \ge k$ are \emph{$k$-trees}.
For $k$-trees of order $n \ge 2k+2$, we characterize all extremal graphs for the upper bound.

\vskip 0.2in \noindent {\emph{keywords}}:
$k$-tree, maximal $k$-degenerate graph, Wiener index   
\end{abstract}

\section{Introduction}
The \emph{Wiener index} of a graph $G$, denoted by $W(G)$,
is the the summation of distances between all unordered vertex pairs of the graph.
The concept was first introduced by Wiener in 1947 for applications in chemistry \cite{W47},
and has been studied in terms of various names and equivalent concepts such as 
the total status \cite{H59}, the total distance  \cite{EJS76},  the transmission \cite{P84},
and the \emph{average distance} (or, \emph{mean distance}) \cite{DG77}.

A graph with a property $\mathcal{P}$ is called \emph{maximal} if it is complete or if adding an edge between any two non-adjacent vertices results 
in a new graph that does not have the property $\mathcal{P}$. 
Finding bounds on Wiener indices of maximal planar graphs of a given order has attracted attention recently, see \cite{CC19, CDOS19+}.
For a maximal planar graph of order $n \ge 3$, its Wiener index has a sharp lower bound $n^2 -4n +6$.
An \emph{Apollonian network} is a chordal maximal planar graph.
Wiener indices of Apollonian networks of order $n \ge 3$  have a sharp upper bound $\lfloor \frac{1}{18} (n^{3} +3n^{2}) \rfloor$,
which also holds for maximal planar graphs of order $3 \le n \le 10$, and was conjectured to be valid for all $n \ge 3$ in \cite{CC19}.
It was shown \cite{CDOS19+} that if $G$ is a $k$-connected maximal planar graph of order $n$,  then the mean distance 
$\mu(G) =\frac{W(G)}{{n \choose 2}} \le \frac{n}{3 k} + O(\sqrt{n})$ for $k \in \{3, 4, 5\}$ and
the coefficient of $n$ is the best possible. 

Let $k$ be a positive integer. 
A graph is $k$-degenerate if its vertices can be successively deleted so that when deleted, they have degree at most $k$.  
Note that Apollonian networks are maximal $3$-degenerate graphs.
In this paper, we provide sharp lower and upper bounds for Wiener indices of maximal $k$-degenerate graphs of order $n$ 
and some extremal graphs for all $n \ge k \ge 1$.
The lower and upper bounds on Wiener indices are equal for maximal $k$-degenerate graphs whose order implies 
that they have diameter at most $2$.
The extremal graphs for the lower bound have a nice description for $2$-trees.
Maximal $k$-degenerate graphs with diameter at least $3$ have order at least $2k+2$.
For $k$-trees of order $n \ge 2k+2$, we charcterize all extremal graphs whose Wiener indices attain the upper bound. 
Our results generalize well-known sharp bounds on Wiener indices of some important classes of graphs 
such as trees and Apollonian networks.

\section{Preliminaries}

All graphs considered in the paper are simple graphs without loops or multiple edges.
Let $G$ be a graph with vertex set $V(G)$ and edge set $E(G)$. 
Then the order of $G$ is $n=|V(G)|$ and the size of $G$ is $|E(G)|$.
Let $K_n$ and $P_n$ denote the clique and the path of order $n$ respectively.
Let $\overline{K}_n$ be the compliment of $K_n$, that is, the graph on $n$ isolated vertices.
Let $G + H$ be the graph obtained from $G$ and $H$ by adding all possible edges between vertices of $G$ and vertices of $H$. 
A complete bipartite graph $K_{r,s}$ is $\overline{K}_r+\overline{K}_s$.

A graph is \emph{connected} if there is a path between any two vertices of the graph.
The  \emph{distance} between two vertices $u,v$ of a graph $G$ is the length of a shortest path 
joining $u$ and $v$ in $G$, and denoted by $d_{G}(u,v)$. 
The distance between two vertices from different components is infinite if $G$ is disconnected.
The \emph{eccentricity} $e_{G}(u)$ of a vertex $u$ in $G$ is the maximum distance between $u$ and other vertices of $G$. 
The set of all vertices with distance 
$i$ from the vertex $u$ in $G$ is denoted by $N_{G}(u,i)$ for $1 \leq i \leq e_{G}(u)$.
In particular, the set of all vertices adjacent to vertex $u$ in $G$ is denoted by $N_G(u)$, 
and its cardinality $|N_G(u)|$ is called the degree of vertex $u$.
The \emph{diameter} of $G$, denoted by $diam(G)$, is  the maximum distance between any two vertices of $G$.
A subgraph $H$ of $G$ is said to be \emph{isometric} in $G$ if $d_{H}(x,y)=d_{G}(x,y)$ for any two vertices $x,y$ of $H$.
The \emph{status} (or, \emph{transmission}) of a vertex   $u$ in $G$, denoted by $\sigma_G(u)$, 
is the summation of the distances between $u$ and all other vertices in $G$. 

\begin{lemma}\cite{BH90, EJS76}\label{L:Preliminary}
Let $G$ be a connected graph. Then

(i) $W(G) \ge 2 {n \choose 2} - |E(G)|$, and the equality holds if and only if $diam(G) \le 2$.

(ii) $W(G) \le W(G-v) + \sigma_G(v)$ for any vertex $v$ of $G$, and the equality holds if and only if $G-v$ is isometric in $G$.

(iii) $W(G)=\sum\limits_{i=1}^{diam(G)} i \cdot d_{i}$,
where $d_{i}$ is the number of unordered vertex pairs with distance $i$ in $G$. 
\end{lemma}

We are interested in $k$-degenerate graphs and maximal $k$-degenerate graphs, introduced in \cite{LW70}. 
A subclass of maximal $k$-degenerate graphs called $k$-trees  \cite{BP69} is particularly important. 
A \emph{$k$-tree} is a generalization for the concept of a tree and can be defined recursively: a clique $K_k$ of order $k \geq 1$ is a $k$-tree,
and any $k$-tree of order $n+1$ can be obtained from a $k$-tree of order $n \geq k$ by adding a new vertex adjacent to all vertices of a clique of order $k$,
which is called the \emph{root} of the newly added vertex, and we say that the newly added vertex is \emph{rooted} at the specific clique.
By definitions, the order of a maximal $k$-degenerate graph can be any positive integer, while the order of a $k$-tree is at least $k$.
A graph  is a \emph{$k$-tree} if and only if it is a chordal maximal $k$-degenerate graph of order $n \ge k$ \cite{B12}. 
A graph is maximal $1$-degenerate if and only if it is a tree \cite{LW70}.
It is known \cite{P86} that $2$-trees form a special subclass of planar graphs extending the concept of maximal outerplanar graphs,
and maximal outerplanar graphs are the only 2-trees that are outerplanar.
Planar 3-trees are just Apollonian networks.

The \emph{$k$-th power of a path} $P_n$, denoted by $P_n^k$, 
has the same vertex set as $P_n$ and two distinct vertices $u$ and $v$ are adjacent in $P_n^k$ 
if and only if their distance in $P_n$ is at most $k$. Note that the order $n$ of $P_n^k$ can be any positive integer.
When $n \ge k$, $P_n^k$ is a special type of $k$-tree.
For $n \ge 2$,  $P_n^k$ is an extremal graph for the upper bound on Wiener indices of maximal $k$-degenerate graphs of order $n$.

A graph is called  \emph{$k$-connected} if the removal of 
any $k-1$ vertices of the graph does not result in a disconnected or trivial graph.
It is well-known that for a $k$-connected graph $G$ of order $n$, $diam\left(G\right)\leq\frac{n-2}{k}+1$.
Since maximal $k$-degenerate graphs of order $n\geq k+1$ are $k$-connected
\cite{LW70}, this bound holds for them, and a characterization
of the extremal graphs (among maximal $k$-degenerate graphs)  appears
in \cite{B12}.

\begin{lemma} \label{L:Upperbound-Status} \cite{CC13, FKM89} 
Let  $G$ be a $k$-connected graph of order $n \ge k+1$  and $k \ge 1$. 
Then $\sigma_G(x) \le  (\lfloor \frac{n - 2}{k} \rfloor +1) (n-1-\frac{k}{2}\lfloor \frac{n-2}{k} \rfloor)$ for any vertex $x$ of $G$.
Moreover, $\sigma_G(x)$
attains the upper bound if and only if $x$ satisfies both properties: 
(i) $e_G(x)=\mbox{diam}(G) = \lfloor \frac{n-2}{k} \rfloor +1$, and 
(ii) $|N_G(x,i)|=k$ for all $1 \le i \le  \lfloor \frac{n-2}{k} \rfloor$.
\end{lemma}

If the graphs in consideration are maximal $k$-degenerate graphs, 
then the upper bound on vertex status in Lemma \ref{L:Upperbound-Status} 
can be achieved by any degree-$k$ vertex of $P_n^k$ for all $n \ge k+1$  and $k \ge 1$.
Furthermore, the extremal graphs are exactly paths $P_n$ when $k=1$.
If $k \ge 2$, then the extremal graphs can be different from $P_n^k$ \cite{B12}.

\section{Sharp Bounds}

\begin{theorem} \label{T:WG)-Lowerbound}
Let $G$ be a $k$-degenerate graph of order $n \ge k \ge 1$.
Then \[W(G) \ge  n^2 - (k+1)n + {k+1 \choose 2}.\] The equality holds if and only if $G$ is maximal $k$-degenerate with $diam(G) \le 2$.
\end{theorem}
\proof By Lemma \ref{L:Preliminary} (i), $W(G) \geq 2 {n \choose 2} - |E(G)|$ and the equality holds if and only if $G$ has diameter at most $2$.
By Proposition 3 in \cite{LW70}, a $k$-degenerate graph $G$ of order $n \geq k$ has $|E(G)| \leq k n - {k+1 \choose 2}$. 
Moreover,  a $k$-degenerate graph $G$ of order $n \geq k$ is maximal if and only if $|E(G)| = k n - {k+1 \choose 2}$, \cite{B12}.  
Therefore, $W(G) \geq n(n-1)-k n + {k+1 \choose 2}  = n^2 - (k+1)n + {k+1 \choose 2}$, and 
the equality holds exactly when $G$ is maximal $k$-degenerate with $diam(G) \leq 2$. \qed\\

This bound is sharp since for $k\leq n\leq k+1$, the only maximal
$k$-degenerate graph is $K_{n}$. For $n\geq k+2$, $ K_k + \overline{K}_{n-k}$ achieves the bound.

\begin{theorem} \label{T:W(G)-Upperbound}
Let $G$ be a maximal $k$-degenerate graph of order $n \geq 2$ and $D=\left\lfloor \frac{n-2}{k}\right\rfloor$. 
Then 
\[
W\left(G\right)\leq W\left(P_{n}^{k}\right)=\sum_{i=0}^{D}\binom{n-ik}{2}=\binom{n}{2}+\binom{n-k}{2}+...+\binom{n-Dk}{2}.
\]
\end{theorem}
\proof We show $W\left(G\right)\leq W\left(P_{n}^{k}\right)$ using induction on order $n$. When $2 \leq n\leq k+2$,
$P_{n}^{k}$ is the only such graph, so it is extremal. Let $G$ be
a maximal $k$-degenerate graph of order $n\geq k+3$, and assume
the result holds for all maximal $k$-degenerate graphs of smaller
orders. By \cite{LW70}, $G$ has a vertex $v$ of degree $k$ and $G-v$ is a maximal $k$-degenerate
graph. Thus $W\left(G-v\right)\leq W\left(P_{n-1}^{k}\right)$.

Label vertices of $P_{n}^{k}$ along the path $P_n$ as $v_1, v_2, \ldots, v_n$ where $n \geq k+3$.
It is clear that $P_n^k$ is $k$-connected and $\sigma_{P_{n}^{k}}\left(v_{n}\right)$  achieves the bound in Lemma \ref{L:Upperbound-Status}. 
By Lemma \ref{L:Preliminary}(iii),
$W\left(G\right)\leq W\left(G-v\right)+\sigma_{G}\left(v\right)
\leq W\left(P_{n}^{k}-v_{n}\right)+\sigma_{P_{n}^{k}}\left(v_{n}\right)=W\left(P_{n}^{k}\right)$.

Note $W\left(P_{n}^{k}\right)=\binom{n}{2}$ when $2\leq n\leq k+1$,
so that the formula holds then. In $P_{n}$, there are $n-i$ pairs
of vertices with distance $i$. Now distances $rk-k+1$ through $rk$
in $P_{n}$ become $r$ in $P_{n}^{k}$. Since $diam\left(P_{n}^{k}\right)=D+1$,
by Lemma \ref{L:Preliminary}(iii),
\begin{align*}
W\left(P_{n}^{k}\right)= & 1\left(n-1\right)+...+1\left(n-k\right)\\
& +2\left(n-k-1\right)+...+2\left(n-2k\right)\\
 & +3\left(n-2k-1\right)+...+3\left(n-3k\right)\\
 & +...\\
 & +D\left(n-\left(D-1\right)k-1\right)+...+D\left(n-Dk\right)\\
 & +\left(D+1\right)\left(n-Dk-1\right)+...+\left(D+1\right)1\\
= & \left(n-1+...+1\right)+\left(n-k-1+...+1\right)+\left(n-2k-1+...+1\right)\\
 & +...+\left(n-\left(D-1\right)k-1+...+1\right)+\left(n-Dk-1+...+1\right)\\
= & \binom{n}{2}+\binom{n-k}{2}+\binom{n-2k}{2}+...+\binom{n-\left(D-1\right)k}{2}+\binom{n-Dk}{2}
\end{align*}
\qed\\

We now provide a closed form expression for $W\left(P_{n}^{k}\right)$
for all $n\geq2$.
\begin{corollary} \label{C:$P_n^k$}
Let $n\geq2$ and $n-2\equiv i\mod k$ for $0\leq i\leq k-1$. Then
\[
W\left(P_{n}^{k}\right)=\frac{n^{3}}{6k}+\frac{\left(k-1\right)n^{2}}{4k}+\frac{\left(k-3\right)n}{12}+\frac{-2i^{3}+3i^{2}\left(k-3\right)-i\left(k^{2}-9k+12\right)-2k^{2}+6k-4}{12k}.
\]
\end{corollary}

\proof We have 

\begin{align*}
W\left(P_{n}^{k}\right)= & \sum_{i=0}^{D}\binom{n-ik}{2}=\sum_{i=0}^{D}\frac{1}{2}\left(n-ik\right)\left(n-ik-1\right)\\
= & \sum_{i=0}^{D}\left[\left(\frac{n^{2}}{2}-\frac{n}{2}\right)+\left(\frac{k}{2}-kn\right)i+\frac{k^{2}}{2}i^{2}\right]\\
= & \sum_{i=0}^{D}\left(\frac{n^{2}}{2}-\frac{n}{2}\right)+\sum_{i=0}^{D}\left(\frac{k}{2}-kn\right)i+\sum_{i=0}^{D}\frac{k^{2}}{2}i^{2}\\
= & \left(D+1\right)\left(\frac{n^{2}}{2}-\frac{n}{2}\right)+\frac{D\left(D+1\right)}{2}\left(\frac{k}{2}-kn\right)+\frac{D\left(D+1\right)\left(2D+1\right)}{6}\frac{k^{2}}{2}\\
= & \frac{k^{2}}{6}D^{3}+\left(\frac{k}{4}+\frac{k^{2}}{4}-\frac{kn}{2}\right)D^{2}+\left(\frac{k}{4}+\frac{k^{2}}{12}-\frac{n}{2}-\frac{kn}{2}+\frac{n^{2}}{2}\right)D-\frac{n}{2}+\frac{n^{2}}{2}
\end{align*}

Since $D=\left\lfloor \frac{n-2}{k}\right\rfloor $, $n-2=Dk+i$ for
$0\leq i\leq k-1$. Substituting $D=\frac{n-2-i}{k}$ into the above
and simplifying, we obtain the formula.
\qed\\

If $1\leq k\leq5$, this formula can be reduced to 
$W\left(P_{n}^{k}\right)=\left\lfloor \frac{2n^{3}+3\left(k-1\right)n^{2}+k\left(k-3\right)n}{12k}\right\rfloor $.
Formulas for small values of $k$ and the beginnings of the resulting
sequences are given in the following table. These sequences occur
(shifted) in OEIS. For $1\leq k\leq3$, they have many different combinatorial
interpretations.
\begin{onehalfspace}
\begin{center}
\begin{tabular}{|c|c|c|c|}
\hline 
$k$ & $W\left(P_{n}^{k}\right)$ & Sequence & OEIS\tabularnewline
\hline 
\hline 
1 & $\frac{n^{3}-n}{6}$ & 0, 1, 4, 10, 20, 35, 56, 84, 120, 165, ... & A000292\tabularnewline
\hline 
2 & $\left\lfloor \frac{n^{3}+1.5n^{2}-n}{12}\right\rfloor $ & 0, 1, 3, 7, 13, 22, 34, 50, 70, 95, ... & A002623\tabularnewline
\hline 
3 & $\left\lfloor \frac{n^{3}+3n^{2}}{18}\right\rfloor $ & 0, 1, 3, 6, 11, 18, 27, 39, 54, 72, ... & A014125\tabularnewline
\hline 
4 & $\left\lfloor \frac{n^{3}+4.5n^{2}+2n}{24}\right\rfloor $ & 0, 1, 3, 6, 10, 16, 24, 34, 46, 61, ... & A122046\tabularnewline
\hline 
5 & $\left\lfloor \frac{n^{3}+6n^{2}+5n}{30}\right\rfloor $ & 0, 1, 3, 6, 10, 15, 22, 31, 42, 55, ... & A122047\tabularnewline
\hline 
\end{tabular}
\par\end{center}
\end{onehalfspace}

 \section{Extremal Graphs}

Any graph of order $n$ and diameter $1$ is a clique and has Wiener index ${n \choose 2}$.
Any maximal $k$-degenerate graph of diameter 1 is $K_{n}$, $2\leq n\leq k+1$,
which is also $P_{n}^{k}$.
Recall that a graph $G$ of order $n$ and diameter $2$ has $W(G)=n(n-1) - |E(G)|$,
and a maximal $k$-degenerate graph $G$ of order $n \ge k$ has $|E(G)| = k n - {k+1 \choose 2}$.  
Then any maximal $k$-degenerate graph of order $n \ge k$ and diameter $2$
has $W(G) = n(n-1)-k n + {k+1 \choose 2}  = {n \choose 2}+ {n-k \choose 2}$.
Therefore, when $k\leq n\leq2k+1$, the lower bound given in Theorem \ref{T:WG)-Lowerbound} and the upper bound given in Theorem \ref{T:W(G)-Upperbound}
are the same, and any maximal $k$-degenerate graph of order $n$ has this value for its Wiener index.

Maximal $1$-degenerate graphs are just trees and so all maximal $1$-degenerate graphs of diameter $2$ are just stars.
For $k \ge 2$, the graphs $ K_k + \overline{K}_{n-k}$ are maximal $k$-degenerate graphs of diameter $2$, but there are others.

We are able to characterize $2$-trees of diameter $2$. 
But the situation becomes complicated as $k$ gets larger.

\begin{figure}[h]
\begin{center}
\includegraphics[width=10cm]{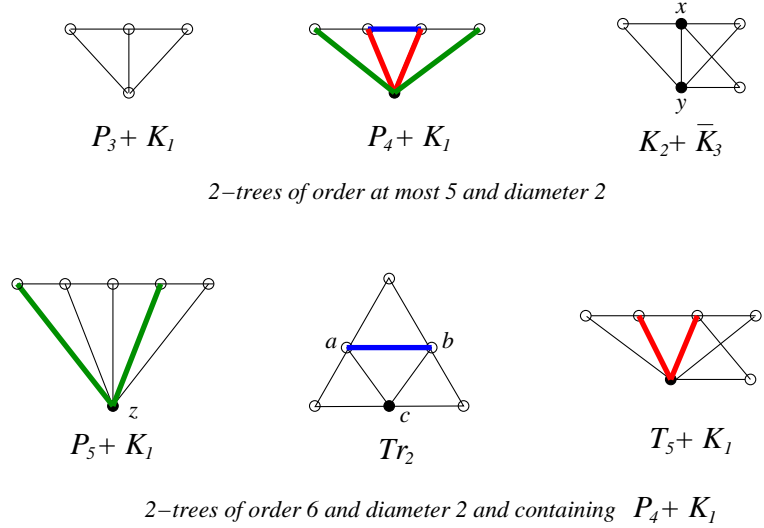}
\caption{Examples of 2-trees.} 
\label{F:2-TreesDiameter2}
\end{center}
\end{figure}

\begin{proposition}\label{P:2-Trees} Let $G$ be a $2$-tree with diameter $2$.
Then $G$ is isomorphic to $T+K_1$ for a tree $T$, 
or a graph formed by adding any number of vertices adjacent to pairs of vertices of $K_3$.
In particular, the maximal outerplanar graphs with diameter $2$ are fans $P_{n-1}+K_1$ and the triangular grid $Tr_2$.
See Figure \ref{F:2-TreesDiameter2}.
\end{proposition}
\proof By its recursive definition, the diameter of $2$-trees cannot decrease as order increases. 
Any $2$-tree with diameter $2$ must have order at least $4$.  
There is a unique $2$-tree with diameter $2$ and order $4$,  $P_4^2=P_3+K_1$.
The $2$-trees of diameter $2$ and order $5$ are $P_5^2 = P_4+K_1$ and $K_2+\overline{K}_3 = K_{1,3}+K_1$.

It is easily seen that $2$-tree not containing $P_4+K_1$ is $K_{1,r}+K_1$ because any additional vertices must be rooted at 
the edge $xy$ of $K_2+\overline{K}_3$, see Figure \ref{F:2-TreesDiameter2}.
Let $G$ be a $2$-tree of order at least 6 and with diameter $2$ containing $P_4+K_1$.  
Then it cannot contain $P_6^2$, the smallest $2$-tree with diameter $3$. 
It is easy to check that $G$ has three possibilities. 

Case 1. $G$ contains $P_5+K_1$. Then any additional vertices must be rooted on edges incident with $K_1$ (the vertex $z$), 
or else it will contain $P_6^2$. 

Case 2.  $G$ contains the triangular grid graph $Tr_2$.
Then the only edges that can be used as roots are those of the central clique $K_3$ (the triangle $abc$), or else it will contain $P_6^2$.

Case 3. $G$ roots all additional vertices on the edges between vertices of degree $3$ and $4$ in $P_4+K_1$. 

Graphs in Case 1 and Case 3 can be described as  $T+K_1$, where $T$ is a tree. 
Graphs in Case 2 are formed by adding vertices rooted at edges from a fixed clique $K_3$.

Maximal planar graphs are exactly the 2-trees that are outerplanar \cite{P86}.  
A graph is outerplanar if and only if it does not contain a subdivision of $K_4$ or $K_{2,3}$ \cite{CH67}.
Thus any maximal outerplanar graph with diameter $2$ is either a fan $P_{n-1}+K_1$ or the triangular grid $Tr_2$.
\qed\\

Since any maximal $k$-degenerate graph of order $n \ge k+1$ is $k$-connected and
$diam(G) \le \lfloor \frac{n-2}{k} \rfloor+1$ for a $k$-connected graph $G$ of order $n$,
any maximal $k$-degenerate graph of diameter at least $3$ has order $n \ge 2k+2$.

\begin{theorem}\label{P:Trees-PathPower} Let $G$ be a $k$-tree of order $n \ge 2k+2$ and $k \ge 1$.
 Then $W(G)=\sum\limits_{i=0}^{\lfloor \frac{n-2}{k} \rfloor} {n-ik \choose 2}$ exactly when $G=P_n^k$.
\end{theorem}
\proof We use induction on order $n$.  
By the recursive definition of a $k$-tree, $G$ can be constructed from a clique $K_{k}$,
and the $i$-th vertex added is adjacent to at least $k-i+1$ vertices of the above clique. 
Thus the smallest order of a $k$-tree with diameter $3$ is $n = 2k + 2$, and the only such $k$-tree is $P_{2k+2}^k$.
So, the result holds for the base case when $n=2k+2$.

Let $G$ be a $k$-tree of order $n\geq2k+3$ that maximizes $W\left(G\right)$,
and assume the result holds for all $k$-trees of order $n-1$.
By the recursive definition of a $k$-tree, $G$ has a vertex $v$ of degree $k$ such that $G-v$ is a $k$-tree. 
By Lemma \ref{L:Preliminary}(ii), 
$W\left(G\right)\leq W\left(G-v\right)+\sigma_{G}\left(v\right)$.

Maximizing $W\left(G-v\right)$ requires that $G-v$ is the extremal graph $P_{n-1}^{k}$. 
Number the vertices of $G-v$ along the path from $1$ to $n-1$.
Since $k$-trees of order at least $k+1$ are $k$-connected, $\sigma_{G}\left(v\right)$ is
maximized when $N_G\left(v\right)=\left\{ 1,2,...,k\right\} $ (or $N_G\left(v\right)=\left\{ n-k,...,n-1\right\} $)
since it achieves the bound in Lemma \ref{L:Upperbound-Status}. When
$n\geq2k+3$, any other choice for $N_G\left(v\right)$ has $\left|N_{G}\left(v,2\right)\right|>k$,
so $\sigma_{G}\left(v\right)$ is not maximized. Thus $G=P_{n}^{k}$, 
and Theorem \ref{T:W(G)-Upperbound} provides the formula.
\qed\\

Note that for $k>1$, there is a unique extremal graph for $k$-trees to achieve the upper bound in Theorem \ref{T:W(G)-Upperbound} 
when $k \leq n\leq k+2$ or $n\geq2k+2$, but not when $k+3\leq n\leq2k+1$.

By Theorem \ref{T:WG)-Lowerbound}, 
Theorem \ref{T:W(G)-Upperbound} and Corollary \ref{C:$P_n^k$}, we have the following 
sharp bounds on Wiener indices of maximal $k$-degenerate graphs for $1 \le k \le 3$.

\begin{corollary} Let $G$ be a maximal $k$-degenerate graph of order $n\geq k \geq 1$. 
\begin{enumerate}
\item If $k=1$,  then $G$ is a tree and 
$n^2 - 2 n  + 1 \leq W(G) \leq \frac{n^3}{6} - \frac{n}{6}.$
The extremal graphs for the bounds are exactly $K_1+\overline{K}_{n-1}$ and $P_n$ respectively, see \cite{EJS76}.

\item If $k=2$, then  
$n^2 - 3 n  +3 \leq W(G) \leq \frac{n^3}{12} + \frac{n^2}{8}- \frac{n}{12} - \frac{1}{16}+\frac{(-1)^n}{16}.$ 

For $2$-trees, the extremal graphs for the lower bound are characterized in Proposition \ref{P:2-Trees};
the extremal graphs for the upper bound are $P_{n}^{2}$ and $K_{2}+\overline{K}_{3}$ (of order 5), 
see Theorem \ref{P:Trees-PathPower}. 

For maximal outerplanar graph of order $n\geq3$ (that is, outerplanar $2$-trees),
the extremal graphs for the lower bound are fans $P_{n-1}+K_1$
and the triangular grid graph $Tr_2$ if $n=6$; and the extremal graphs for the upper bound are $P_n^2$.

\item If $k=3$, then  
$n^2 - 4 n +6 \leq W(G) \leq \lfloor\frac{n^3}{18} + \frac{n^{2}}{6}\rfloor$.

For $3$-trees, it is easily checked that the extremal
graphs for the upper bound are $P_{n}^{3}$, $K_{3}+\overline{K}_{3}$ of order 6 and four others
of order $7$ which are $K_{3}+\overline{K}_{4}$, $K_{2}+T_5$, where
$T_5$ is the tree of order 5 that is neither a path nor a star, $P_{5}+K_{2}$,
and the graph formed from $K_{4}$ by adding degree $3$ vertices inside $3$ regions. 
See Figure \ref{F:3-TreesDiameter2}.

For Apollonian networks (planar $3$-trees), the upper bound was given in \cite{CC19}.  
The extremal graphs for the upper bound are $P_{n}^{3}$
and the last two graphs of order $7$ in Figure \ref{F:3-TreesDiameter2}.
\end{enumerate}
\begin{figure}[h]
\begin{center}
\includegraphics[width=10cm]{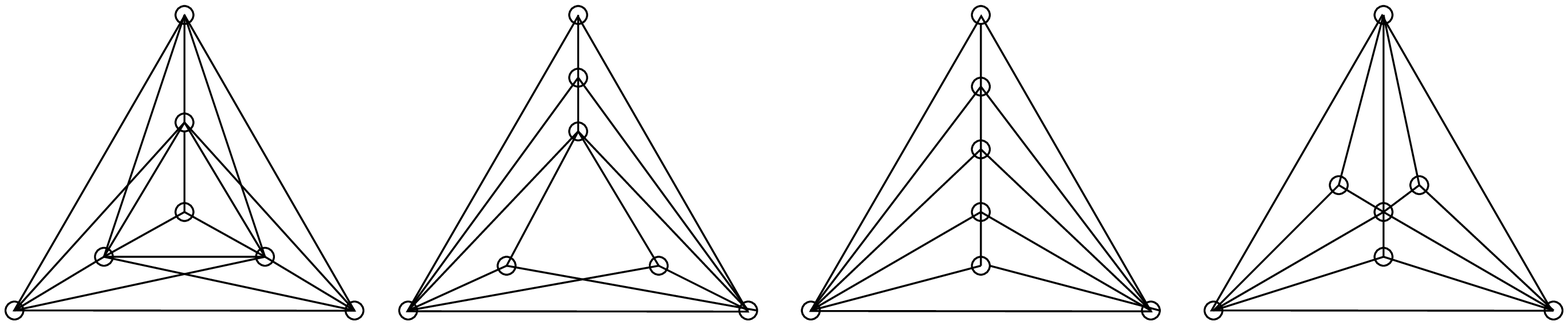}
\caption{Examples of 3-trees of order $7$.} 
\label{F:3-TreesDiameter2}
\end{center}
\end{figure}
\end{corollary}



\end{document}